\documentclass[12pt]{article}

\usepackage{amssymb}
\usepackage{amsmath}
\usepackage{lscape}
\usepackage{url}
\usepackage{fullpage}

\input xy
\xyoption{all}

\makeatletter
\newcommand\Small{\@setfontsize\small\@ixpt{11}}   
\makeatother

\newcommand{\matX}{\mathbf{X}}
\newcommand{\matI}{\mathbf{I}}
\newcommand{\matM}{\mathbf{M}}
\newcommand{\II}{\mathcal{I}}
\newcommand{\XX}{\mathcal{X}}

\newcommand{\MM}{\mathcal{M}}
\newcommand{\TT}{\mathcal{T}}
\newcommand{\jalg}{\mathbf{H}_3(\OO)}

\newcommand{\Aut}{\textrm{Aut}}
\renewcommand{\Im}{\textrm{Im}}
\newcommand{\stab}{\textrm{stab}}

\newcommand{\dotBtz}{\dot B^2_{tz} - \dot B^3_{tz}}

\newcommand{\KK}{\mathbb{K}}
\newcommand{\RR}{\mathbb{R}}
\newcommand{\CC}{\mathbb{C}}
\newcommand{\HH}{\mathbb{H}}
\newcommand{\OO}{\mathbb{O}}

\newcommand{\ga}{\mathfrak{g}}

\renewcommand{\aa}{\mathfrak{a}}
\newcommand{\bb}{\mathfrak{b}}
\newcommand{\cc}{\mathfrak{c}}
\newcommand{\dd}{\mathfrak{d}}
\newcommand{\ee}{\mathfrak{e}}
\newcommand{\ff}{\mathfrak{f}}
\newcommand{\gtwo}{\mathfrak{g}_2}

\newcommand{\sla}{\mathfrak{sl}}
\newcommand{\so}{\mathfrak{so}}
\newcommand{\su}{\mathfrak{su}}
\renewcommand{\sp}{\mathfrak{sp}}

\newcommand{\uonem}{\mathfrak{u}(-1)}

\newcommand{\SL}{SL}
\newcommand{\SO}{SO}
\newcommand{\SU}{SU}
\newcommand{\Uone}{U(1)}
\newcommand{\EE}{E}
\newcommand{\FF}{F}
\newcommand{\Gtwo}{G_2}

\newcommand{\suC}{\su(3)_\CC}
\newcommand{\suH}{\su(2)_\HH}
\newcommand{\soH}{\so(3)_\HH}
\newcommand{\SUC}{\SU(3)_\CC}
\newcommand{\SUH}{\SU(2)_\HH}
\newcommand{\SOH}{\SO(3)_\HH}

\begin{document}


\title{\boldmath\textbf{$E_6$, the Group:\\
	The structure of $\SL(3,\OO)$}}

\author{%
Aaron Wangberg\\
Department of Mathematics \& Statistics\\
Winona State University\\
Winona, MN 55987\\
\texttt{awangberg@winona.edu}
\and
Tevian Dray\\
Department of Mathematics\\
Oregon State University\\
Corvallis, OR  97331\\
\texttt{tevian@math.oregonstate.edu}
}

\date{\today}

\maketitle

\begin{abstract}
We present the subalgebra structure of $\sla(3,\OO)$, a particular real form of
$e_6$ chosen for its relevance to particle physics and its close relation to
generalized Lorentz groups.  We use an explicit representation of the Lie
group $\SL(3,\OO)$ to construct the multiplication table of the corresponding
Lie algebra $\sla(3,\OO)$.  Both the multiplication table and the group are then
utilized to find various nested chains of subalgebras of $\sla(3,\OO)$, in which
the corresponding Cartan subalgebras are also nested where possible.  Because
our construction involves the Lie group, we simultaneously obtain an explicit
representation of the corresponding nested chains of subgroups of $\SL(3,\OO)$.
\end{abstract}

\section{Introduction}

The group $E_6$ has a long history of applications in
physics~\cite{Jordan,JNW,Albert}, and is a candidate gauge group for a Grand
Unified Theory~\cite{Georgi}.  A description of the group $E_{6(-26)}$ as
$\SL(3,\OO)$ was given in~\cite{Denver}, generalizing the interpretation of
$\SL(2,\OO)$ as (the double cover of) $\SO(9,1)$ discussed in~\cite{Lorentz}.
An interpretation combining spinor and vector representations of the Lorentz
group in 10 spacetime dimensions was described in~\cite{York}.  In this paper,
we fill in some further details of the structure of $\SL(3,\OO)$, in the
process obtaining nested chains of subgroups that respect this Lorentzian
structure.

We begin by reviewing the construction of both $\SL(2,\OO)$ and $\SL(3,\OO)$ at
the group level in Section~\ref{Group}, then describe the construction of the
Lie algebra $\sla(3,\OO)$ in Section~\ref{Algebra}.  In Section~\ref{Chains}, we
use this information to construct various chains of subgroups and subalgebras,
some but not all of which are simple, and in Section~\ref{conclusion} we
discuss our results.

\section{The Group}
\label{Group}

\subsection{Lorentz transformations}
\label{Lorentz}

\begin{table}[btp]
\begin{center}
$$
\begin{array}[t]{|l|c|}
\hline
\hbox{Boosts} &
\begin{array}{ccl}
\noalign{\vspace{0.05in}}
B_{tz} & t \longleftrightarrow z &
 \matM = \left( \begin{array}{cc}
	\exp\left(\frac{\alpha}{2}\right) & 0 \\
	0 & \exp\left(-\frac{\alpha}{2}\right) \\
	\end{array} \right) \\
\noalign{\vspace{0.05in}}
B_{tx} & t \longleftrightarrow x &
 \matM = \left( \begin{array}{cc}
    \cosh\left(\frac{\alpha}{2}\right) & \sinh\left(\frac{\alpha}{2}\right) \\
    \sinh\left(\frac{\alpha}{2}\right) & \cosh\left(\frac{\alpha}{2}\right)
    \end{array} \right) \\
\noalign{\vspace{0.05in}}
B_{tq} & t \longleftrightarrow q &
 \matM = \left( \begin{array}{cc}
    \cosh\left(\frac{\alpha}{2}\right) & -q \sinh\left(\frac{\alpha}{2}\right)\\
    q \sinh\left(\frac{\alpha}{2}\right) & \cosh\left(\frac{\alpha}{2}\right) 
    \end{array} \right) \\
\noalign{\vspace{0.05in}}
\end{array} \\
\hline
\hline
\begin{minipage}{0.85in}
Simple \\ Rotations
\end{minipage} &
\begin{array}{ccl}
\noalign{\vspace{0.05in}}
R_{xq} & x \longleftrightarrow q &
 \matM = \left( \begin{array}{cc}
	\exp\left(-\frac{q\alpha}{2}\right) & 0 \\
	0 & \exp\left(\frac{q\alpha}{2}\right)
	\end{array} \right) \\
\noalign{\vspace{0.05in}}
R_{xz} & x \longleftrightarrow z &
 \matM = \left( \begin{array}{cc}
    \cos\left(\frac{\alpha}{2}\right) & \sin\left(\frac{\alpha}{2}\right) \\
   -\sin\left(\frac{\alpha}{2}\right) & \cos\left(\frac{\alpha}{2} \right) 
    \end{array} \right) \\
\noalign{\vspace{0.05in}}
R_{zq} & q \longleftrightarrow z &
 \matM = \left( \begin{array}{cc}
    \cos\left(\frac{\alpha}{2}\right) & q \sin\left(\frac{\alpha}{2}\right) \\
    q \sin\left(\frac{\alpha}{2}\right) & \cos\left(\frac{\alpha}{2}\right) 
    \end{array} \right) \\
\noalign{\vspace{0.05in}}
\end{array} \\
\hline
\hline
\begin{minipage}[c]{0.85in}
\vspace{0.1in}
Transverse \\ Rotations
\vspace{0.15in}
\end{minipage} & 
\begin{array}{ccl}
R_{p,q} & p \longleftrightarrow q &
 \matM_1 = -p \; \matI_2 \\
  & &\matM_2 = \left(\; \cos \left(\frac{\alpha}{2}\right) p
	+ \sin \left( \frac{\alpha}{2} \right) q \; \right) \matI_2
\end{array} \\
\hline
\end{array}
$$
\caption{Finite octonionic Lorentz transformations.  The group transformation
is given by \hbox{$\matX\longmapsto\matM\matX\matM^\dagger$} for boosts and
simple rotations, and by
$\matX\longmapsto\matM_2(\matM_1\matX\matM_1^\dagger)\matM_2^\dagger$ for
transverse rotations.  The parameters $p$ and $q$ are imaginary unit
octonions.}
\label{Finite}
\end{center}
\end{table}

A $2\times2$ Hermitian matrix
\begin{equation}
\matX = \left( \begin{array}{cc} t+z& x-q\\ x+q& t-z\end{array} \right)
\end{equation}
with $t,z,x\in\RR$ and pure imaginary $q\in\KK=\RR,\CC,\HH,\OO$, is a
representation of an ($m+1$)-dimensional spacetime vector for
$m+1=||\KK||+2\in\{3,4,6,10\}$.  In this setting, the squared Lorentzian norm
of $\matX$ is given by $\det\matX$.  Lorentz transformations preserve
$\det\matX$ and must also preserve the Hermiticity of $\matX$.  Any Lorentz
transformation can be described as the composition of maps of the form
\begin{equation}
\matX \mapsto \matM \matX \matM^\dagger
\label{Mtwo}
\end{equation}
for certain \textit{generators} $\matM$.  In the octonionic case, the
determinant-preserving transformations of the form~(\ref{Mtwo}) constitute
$\SL(2,\OO)$, the (double cover of the) Lorentz group $\SO(9,1)$.  We adopt the
explicit set of generators constructed by Manogue and Schray~\cite{Lorentz},
as given in Table~\ref{Finite}.

An important feature of these transformations is that the \textit{transverse
rotations} between octonionic units require \textit{nesting}; the lack of
associativity prevents one from combining the given transformations of the
form~(\ref{Mtwo}) into a single such transformation.  We will return to this
point in Section~\ref{Algebra}.

The exceptional Jordan algebra~$\jalg$, also known as the Albert algebra,
consists of~\hbox{$3\times3$} octonionic Hermitian matrices under the Jordan
product, and forms a~$27$-dimensional representation of~$E_6$~\cite{corrigan},
which is precisely the group that preserves the determinant of Jordan
matrices; in this sense, $E_6=\SL(3,\OO)$.  There are three natural ways to
embed a $2\times2$ Hermitian matrix in a $3\times3$ Hermitian matrix, as
illustrated in Table~\ref{Types}, which we refer to as \textit{types}.
Furthermore, $\SL(2,\OO)$ sits inside $\SL(3,\OO)$ under the identification
\begin{equation}
\matM \longmapsto
\MM = \left( \begin{array}{c|c}
	\matM & 0 \\ \hline 0 & 1
	\end{array}\right)
\label{typeI}
\end{equation}
If we take $\XX\in\jalg$ to be of type~1, as per Table~\ref{Types}, then
under the transformation
\begin{equation}
\XX \mapsto \MM \XX \MM^\dagger
\label{Mthree}
\end{equation}
we recover not only the \textit{vector} transformation~(\ref{Mtwo}) on
$\matX$, but also the \textit{spinor} transformation
\(
\theta \mapsto \MM\theta
\)
on the 2-component octonionic column $\theta$.%
\footnote{Further discussion of ``vectors'' and ``spinors'' can be found
in~\cite{York}.}

\begin{table}[tbp]
\begin{center}
\begin{tabular}{ccccc}
\textrm{Type 1} & & \textrm{Type 2} & & \textrm{Type 3}\\
$\left( \begin{array}{c|c}
	\matX & \theta \\ \hline \theta^\dagger & \cdot
	\end{array}\right)$ & &
$\left( \begin{array}{c|c}
	\cdot & \theta^\dagger \\ \hline \theta & \matX
	\end{array}\right)$ & &
$\left( \begin{array}{c|c|c}
	\matX_{2,2} & \theta_2 & \matX_{2,1} \\
	\hline \overline{\theta_2} & \cdot & \overline{\theta_1} \\
	\hline \matX_{1,2} & \theta_{1} & \matX_{1,1} \end{array}\right)$ \\
\end{tabular}
\caption{Three natural locations of a vector~$\matX$, a spinor~$\theta$, and a
dual spinor~$\theta^\dagger$ in~$\XX\in E_6$.}
\label{Types}
\end{center}
\end{table}

We generalize this construction to all three types.  We write $M^1$ (instead
of $\MM$) for the type~1 version of $\MM$, as defined by~(\ref{Mthree}).  Then
type~2 and~3 versions of $\MM$ can be obtained as
\begin{equation}
M^2 = \TT M^1 \TT^\dagger \hspace{1.5cm}
M^3 = \TT M^2 \TT^\dagger \hspace{1.5cm}
\label{typetran}
\end{equation}
so that the group transformation
\begin{equation}
\TT
  = \left( \begin{array}{ccc}
	0 & 0 & 1 \\ 1 & 0 & 0 \\ 0 & 1 & 0
	\end{array} \right)
  \in E_6
\end{equation}
cyclically permutes the 3 types.  We discuss \textit{type transformations} of
the form~(\ref{typetran}) in more detail in Section~\ref{Type}

\subsection{A new basis for transverse rotations}

As outlined in~\cite{Denver,York}, $E_6$ can be viewed as the appropriate
union of these 3 copies of \hbox{$\SO(9,1,\RR) = \SL(2,\OO)$}.  But we have
$3\times45=135$ elements, and we need to find a way to reduce this number to
$|E_6|=78$.
We start by constructing a new basis for the transverse rotations in
$\SL(2,\OO)$.

Each transverse rotation $R_{p,q}$ listed in Table~\ref{Finite} rotates a
single plane spanned by the orthogonal imaginary octonions~$p$ and~$q$, and
rotations of the 21 independent planes generate $\SO(7)$.  Since
\hbox{$G_2\subset \SO(7)$}, we choose a basis for $G_2$ and extend it to
$\SO(7)$.  For each basis octonion, say $q=i$, there are three pairs of basis
octonions, in this case $\lbrace j,k \rbrace$, $\lbrace k\ell,j\ell \rbrace$,
$\lbrace \ell,i\ell \rbrace$, which generate quaternionic subalgebras
containing $q$.  We have chosen the ordering of the pairs so that adding $q$
leads to a right-handed, three-dimensional coordinate frame, and so that
$\ell$ only appears (if at all) in the last pair.  A choice of pairs for each
basis octonion that satisfies these conditions is given in Table~\ref{AGS}.
We now define the combinations
\begin{align}
A_i(\alpha) &= R_{j,k}(\alpha) \circ R_{k\ell,j\ell}(-\alpha) \nonumber\\
G_i(\alpha) &= R_{j,k}(\alpha) \circ R_{k\ell,j\ell}(\alpha)
			\circ R_{\ell,i\ell}(-2\alpha) \label{AGSdef}\\
S_i(\alpha) &= R_{j,k}(\alpha) \circ R_{k\ell,j\ell}(\alpha)
			\circ R_{\ell,i\ell}(\alpha) \nonumber
\end{align}
and use the conventions in Table~\ref{AGS} to similarly define $A_q$, $G_q$,
and $S_q$ for the remaining basis octonions.

As we discuss in more detail below, the 14 transformation of the form $A_q$
and $G_q$ generate the group~$G_2$, the 7 transformations $A_q$ together with
$G_\ell$ generate the subgroup $\SU(3)\subset G_2$ which fixes~$\ell$, and all
21 of these transformations, which generate $\SO(7)$, are orthogonal (but not
normalized) at the Lie algebra level.  We will use these properties to
eliminate redundant group generators.

\begin{table}[tbp]
\begin{center}
\begin{tabular}{|c|c|c|c|}
\hline
$q$ & First pair & Second pair & Third pair \\
\hline
\hline
$i$  & $(j , k )$ & $(k\ell , j\ell)$ & $(\ell , i\ell)$ \\
$j$  & $(k , i )$ & $(i\ell , k\ell)$ & $(\ell , j\ell)$ \\
$k$  & $(i , j )$ & $(j\ell , i\ell)$ & $(\ell , k\ell)$ \\
$k\ell$ & $(j\ell , i)$ & $(j , i\ell) $ & $(k , \ell) $ \\
$j\ell$ & $(i , k\ell)$ & $(i\ell , k) $ & $(j , \ell) $ \\
$i\ell$ & $(k\ell , j)$ & $(k , j\ell) $ & $(i , \ell) $ \\
$\ell$  & $(i\ell , i)$ & $(j\ell , j) $ & $(k\ell , k)$ \\
\hline
\end{tabular}
\caption{Quaternionic subalgebras chosen for~$A_q$,~$G_q$, and~$S_q$.}
\label{AGS}
\end{center}
\end{table}

\section{The Lie Algebra}
\label{Algebra}

\subsection{Constructing the algebra}
\label{E6alg}

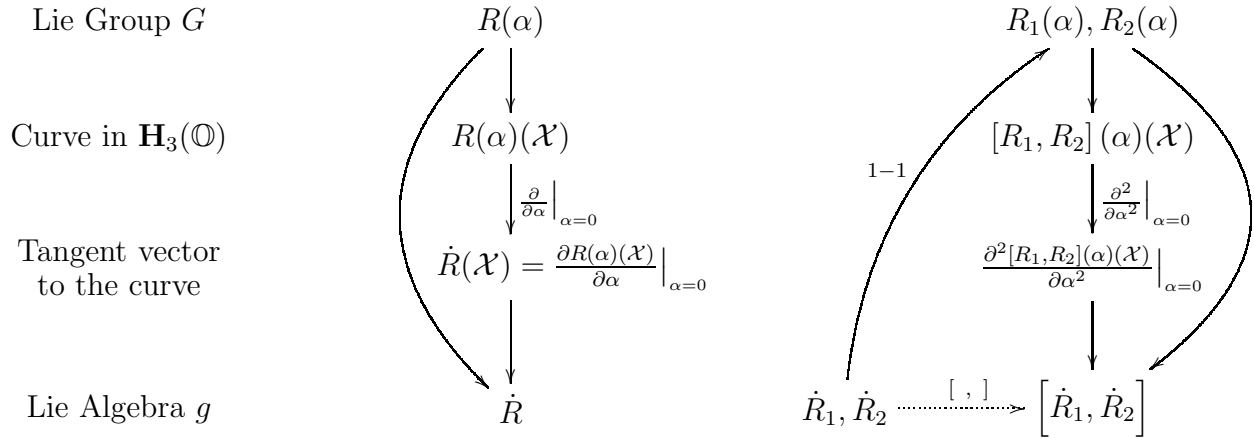
\begin{figure}[tbp]
\begin{center}
\xymatrix@M=2pt{
*++\txt{Lie Group~$G$} &
*++\txt{$R(\alpha)$}\ar@/_3.5pc/[ddd] \ar[d] &
\hspace{1.5cm} &
*++\txt{$ R_1(\alpha),R_2(\alpha) $} \ar[d] \ar@/^5.0pc/[ddd]\\
*++\txt{$\textrm{Curve in } \jalg$} &
*++\txt{$R(\alpha)(\XX)$} \ar[d]^{\frac{\partial }{\partial \alpha} \big|_{_{\alpha = 0}}}
& & *++\txt{$\left[R_1,R_2\right](\alpha)(\XX)$}
	\ar[d]^{\frac{\partial^2 }{\partial \alpha^2}\big|_{_{\alpha = 0}}} \\
*++\txt{$\textrm{Tangent vector}$\\$\textrm{to the curve}$} &
*++\txt{$\qquad\qquad\dot{R}(\XX)=\frac{\partial R(\alpha)(\XX)}
		{\partial \alpha}\big|_{_{\alpha = 0}}$}\ar[d]
& & *++\txt{$\frac{\partial^2 \left[ R_1, R_2 \right](\alpha)(\XX)}
	{\partial \alpha^2}\big|_{_{\alpha = 0}}$} \ar[d] \\
*++\txt{Lie Algebra~$g$} &
*++\txt{$\dot R$} &
*++\txt{$\dot R_1, \dot R_2 $} \ar@/^2pc/[uuur]^{1-1}
	\ar@{.>}[r]^{\left[ \hspace{.15cm}, \hspace{.15cm} \right]} &
*++\txt{$\left[ \dot R_1, \dot R_2 \right]$}\\
}
\caption{Calculating Lie algebra elements and their commutators.}
\label{Commutators}
\end{center}
\end{figure}

We begin by associating each transformation in the Lie group with a vector in
the Lie algebra.  Each of the~$135$ transformations is a one-parameter curve
in the group.  Given a one-parameter curve~$R(\alpha)$ in a classical Lie
group, the traditional method for associating it with the Lie algebra
generator~$\dot R$ is to find its tangent vector
$\dot R = \frac{\partial R(\alpha)}{\partial\alpha} \big|_{_{\alpha = 0}}$
at the identity element in the group.  However, the transverse rotations are
\textit{nested}, that is, they involve more than one operation, and the lack
of associativity prevents one from working with the group elements by
themselves.  Instead, we let our one-parameter transformations $R(\alpha)$ act
on elements $\XX\in\jalg$, producing a curve~$R(\alpha)(\XX)$ in~$\jalg$.  We
then define the Lie algebra element $\dot{R}\in\ee_6$ to be the map taking
$\XX$ to the tangent vector at the identity to this curve in~$\jalg$.  That
is, we have the association indicated in Figure~\ref{Commutators} between the
group transformations and the tangent vectors.

We also use group orbits to construct the commutator of two tangent vectors.
In the traditional approach to the classical matrix groups, the commutator of
the tangent vectors~$\dot R_1$ and~$\dot R_2$ is defined as \hbox{$[ \dot R_1,
\dot R_2 ] = \dot R_1 \dot R_2 - \dot R_2 \dot R_1$}.  However, we are working
in~$\jalg$, not~$\ee_6$.  To find the commutator of the Lie algebra elements
$\dot{R_1}$ and $\dot{R_2}$ associated with curves~$R_1(\alpha)(\XX)$
and~$R_2(\alpha)(\XX)$, we create a new curve in~$\jalg$ defined by
\begin{equation}
[ R_1, R_2](\alpha)(\XX)
  = R_1(\frac{\alpha}{2}) \circ R_2(\frac{\alpha}{2})
	\circ R_1(-\frac{\alpha}{2}) \circ R_2(-\frac{\alpha}{2})(\XX)
\end{equation}
where~$\circ$ denotes composition.  This new path is not a one-parameter
curve, and its first derivative is identically zero at~$\alpha = 0$, but its
second derivative is tangent to the curve $[R_1, R_2](\alpha)(\XX)$
at~$\alpha=0$.  Therefore, we define the commutator of~$\dot R_1$ and~$\dot
R_2$ by the following action on $\jalg$
\begin{equation}
\left[ \dot R_1, \dot R_2 \right] (\XX)
  = \frac12 \frac{\partial^2}{\partial \alpha^2}
	\bigl[R_1, R_2\bigr](\alpha)(\XX) \big|_{_{\alpha = 0}}
\end{equation}
which agrees with the usual definition for matrix Lie groups~\cite{gilmore}.
Our construction of the commutator is summarized in Figure~\ref{Commutators}.

Since we are using the local action of~$\SL(3,\OO)$ on~$\jalg$ to give a
homomorphic image of~$\sla(3,\OO)$, our construction does not lead to a readily
available exponential map giving the group element corresponding
to~$[\dot R_1,\dot R_2]$.  In particular, we are not always able to find the
one-parameter curve whose tangent vector is~$[\dot R_1,\dot R_2]$.

\subsection{Linear dependencies}
\label{Depend}

We shall now give the dependencies among the group transformations by using
linear dependencies among the Lie algebra elements.  In doing so, we will
indicate which transformations can be eliminated, leaving our preferred basis
for the group~$\SL(3,\OO)$ and the algebra~$\sla(3,\OO)$.  Since we are using a
homomorphic image of the Lie algebra~$\sla(3,\OO)$, we check that the indicated
dependencies actually do provide dependencies among the group transformations.

We begin with the transverse rotations.  Among the~$21$
transformations~$A_q$, $G_q$, and~$S_q$ of each type, direct computation
shows that
\begin{equation}
\dot A^1_q
  = \dot A^2_q
  = \dot A^3_q \hspace{2cm}
\dot G^1_q
  = \dot G^2_q
  = \dot G^3_q
\end{equation}
for each basis octonion~$q$.  That is, the transformations~$A_q$ and~$G_q$ are
type independent, allowing us to drop the type designation and simply
write~$\dot A_q$ and~$\dot G_q$.  These fourteen transformations generate
$G_2=\Aut(\OO)$, which is the smallest of the exceptional Lie groups.  We
refer to the type independence of these transformations as \textit{strong
triality}.
When added to the fourteen~$G_2$ transformations, the seven
transformations~$S^a_q$ produce a basis for the~$\SO(7)$ of type~$a$, with
\hbox{$a=1,2,3$}.  However, the transformations~$S^a_q$ are not independent,
since
\begin{equation}
\dot S^1_q + \dot S^2_q + \dot S^3_q = 0
\end{equation}
Hence, the union of any two of the~$\SO(7)$ subgroups contains the third.  In
particular, we may use the group transformations generated by~$S^a_q$ of
type~$1$ and type~$2$ to generate the type~$3$ transformations generated
by~$S^3_q$.  These linear dependences have reduced our~$3\times21=63$
transverse rotations by~$28+7=35$, trimming our original~$135$ transformations
down to~$100$.

Turning to $\SO(8)$, we have the relations
\begin{align}
0 &= \dot R^1_{xq} + \dot R^2_{xq} + \dot R^3_{xq} \nonumber\\
\dot R^2_{xq} &= -\frac{1}{2}\dot R^1_{xq} - \frac{1}{2}\dot S^1_q 
\label{Srels}\\
\dot S^2_q &= \frac{3}{2} \dot R^1_{xq} - \frac{1}{2}\dot S^1_q \nonumber
\end{align}
which allow us to eliminate a further 21 transformations.  We have in fact
expressed all $\SO(8)$ transformations of types~2 and~3 in terms of $\SO(8)$
transformations of type~1; in this sense, there is only one~$\SO(8)$!  Again,
this is a result of \textit{triality}.

Having reduced the~$135$ transformations to~$100$ and then by another~$21$
to~$79$, we are left with~$52$ rotations, which preserve the trace of
\hbox{$\XX\in\jalg$}, and which form the Lie group~$F_4=\SU(3,\OO)$.
Among the remaining~$27$ boosts, we expect only one additional linear
dependency, which turns out to be
\begin{equation}
\dot B^1_{tz} + \dot B^2_{tz} + \dot B^3_{tz} = 0
\end{equation}
which we use to eliminate~$\dot B^2_{tz}$ and~$\dot B^3_{tz}$ in favor of the
combination~$\dotBtz$.  The resulting $78$ Lie algebra elements are indeed
independent, and turn out to be orthogonal (but not normalized) with respect
to the Killing form.

We have therefore constructed both the group $E_6=\SL(3,\OO)$, and its Lie
algebra \hbox{$\ee_6=\sla(3,\OO)$}; the complete commutation table
for~$\sla(3,\OO)$ can be found online at~\cite{commutation_table_online}.  In
retrospect, the counting is easy: There is one $\SO(8)$ (28 elements), 3 types
of each of the remaining elements of $\SO(9)$ (24 elements, yielding $F_4$),
and 3 types of the 9 boosts, with one final dependency, yielding 26 boosts in
all.

Our basis can be simplified slightly by noticing that~(\ref{Srels}) implies
\begin{equation}
\dot S^1_q = \dot R^3_q - \dot R^2_q
\end{equation}
where the operations on the RHS commute.  Thus, the diagonal \textit{phase}
transformations $S^1_q$ can in fact be constructed \textit{without} nesting,
which however is essential for the $G_2$ transformations $A_q$ and $G_q$.
This provides another way to count the basis of $\ee_6$: There are 64
independent trace-free $3\times3$ octonionic matrices, $24+14=38$ of which are
anti-Hermitian (infinitesimal rotations), and $24+2=26$ of which are Hermitian
(boosts), together with the 14 nested transformations making up $\gtwo$, for a
total of 78 independent elements in $\ee_6$~\cite{corrigan}.

We can further identify the 6 elements
\begin{equation}
C = \lbrace \dot B^1_{tz}, \dot B^2_{tz}-\dot B^3_{tz}, \dot R^{1}_{x\ell},
	\dot A_\ell, \dot G_\ell, \dot S^1_{\ell} \rbrace
\end{equation}
as a commuting set, and therefore a preferred (orthogonal) basis for the
Cartan subalgebra~$h$.  We call these basis elements the \textit{Cartan
elements} of~$\ee_6$.

\section{Subalgebra Chains}
\label{Chains}

\subsection{Basic subalgebra chains}
\label{LorentzSub}

We begin with a discussion of $\gtwo\subset \so(7)$.  Our basis selects a
preferred $\su(3)$ subalgebra of $\gtwo$, namely the $\gtwo$ transformations
which fix the preferred complex subalgebra of $\OO$ generated by $\ell$.
G\"unaydin denotes the corresponding~$\SU(3)$ subgroup of~$G_2$
as~$\SU(3)^C$~\cite{gunaydin_1974}; we prefer to use the name $\suC$ for this
subalgebra.  Explicitly, we have
\begin{equation}
\suC = \langle \dot A_i, \cdots, \dot A_{i\ell}, \dot A_\ell, \dot G_\ell \rangle
\end{equation}
which is also a subalgebra of the (type~1, say) $\so(6)\subset\so(7)$ that
fixes $\ell$.

Through a conventional choice of $A_\ell$, our basis also selects a preferred
quaternionic subalgebra of $\OO$, generated by $\lbrace k,k\ell,\ell \rbrace$,
and a preferred subalgebra $\suH\subset\suC$ that fixes this quaternionic
subalgebra, namely
\begin{equation}
\suH = \langle \dot A_k, \dot A_{k\ell}, \dot A_\ell \rangle
\end{equation}
Extending to $\so(7)$, there is clearly an $\so(4)$ that fixes $\HH$; we have
\begin{equation}
\so(4) = \so(3) \oplus \so(3) = \suH \oplus
	\langle \dot G_k - \dot S^1_k, \dot G_{k\ell} - \dot S^1_{k\ell},
	\dot G_\ell - \dot S^1_\ell \rangle
\end{equation}
as can be seen by studying Table~\ref{AGS}.  Another interesting $\so(3)$
subalgebra of $\so(7)$ is the complement of this $\so(4)$, an orthogonal basis
for which is given by the combinations $\dot G_q + 2\dot S^1_q$ for
$q\in\Im\HH$.

\begin{figure}[tbp]
{\Small
\begin{center}
\begin{minipage}{6in}
\begin{center}
\[
\xymatrix{
 & *++\txt{$\SL(3,\OO)$\\ $[\ee_6]$} & \\
*++\txt{$\SO(9,1) = \SL(2,\OO)$\\
  $[\dd_5]$}\ar[ur]^{B^2_{tz}} & *++\txt{$\SO(9,1) = \SL(2,\OO)$\\
  $[\dd_5]$}\ar[u]_{B^1_{tz}} & *++\txt{$\SO(9,1) = \SL(2,\OO)$\\
  $[\dd_5]$}\ar[ul]^{B^3_{tz} \to B^1_{tz}}_{B^2_{tz}} \\
*++\txt{$\SO(9) = \SU(2,\OO)$\\
  $[\bb_4]$}\ar[u]^{B^1_{tz}} & *++\txt{$\SO(9) = \SU(2,\OO)$\\
  $[\bb_4]$}\ar[u]^{B^2_{tz}} & *++\txt{$\SO(9) = \SU(2,\OO)$\\
  $[\bb_4]$}\ar[u]^{B^3_{tz}} \\
 & *++\txt{$\SO(8)$\\$[\dd_4]$}\ar[u] \ar[ur] \ar[ul] \\
*++\txt{$\SO(7) = \SU(1,\OO)$\\
  $[\bb_3]$}\ar[ur]^{R^1_{x\ell}} & *++\txt{$\SO(7) = \SU(1,\OO)$\\
  $[\bb_3]$}\ar[u]^{S^2_\ell \to S^1_\ell}_{R^1_{x\ell}} &
    *++\txt{$\SO(7) = \SU(1,\OO)$\\
  $[\bb_3]$}\ar[ul]^{S^3_\ell \to S^1_\ell}_{R^1_{x\ell}} \\
*++\txt{$\SO(6)$\\
  $[\dd_3]$}\ar[u] & *++\txt{$\SO(6)$\\
  $[\dd_3]$}\ar[u] & *++\txt{$\SO(6)$\\
  $[\dd_3]$}\ar[u] \\
*++\txt{$\SO(5)$\\
  $[\bb_2]$}\ar[u]_{G_\ell-S^1_\ell} & *++\txt{$\SO(5)$\\
  $[\bb_2]$}\ar[u]_{G_\ell-S^2_\ell} & *++\txt{$\SO(5)$\\
  $[\bb_2]$}\ar[u]_{G_\ell-S^3_\ell} \\
*++\txt{$\SO(4) = \SO(3)\times\SO(3)$\\
  $[\dd_2=\aa_1\oplus\aa_1]$}\ar[u] & *++\txt{$\SO(4) = \SO(3)\times\SO(3)$\\
  $[\dd_2=\aa_1\oplus\aa_1]$}\ar[u] & *++\txt{$\SO(4) = \SO(3)\times\SO(3)$\\
  $[\dd_2=\aa_1\oplus\aa_1]$}\ar[u] \\
 & *++\txt{$\SO(3)$\\
  $[\bb_1]$}\ar[u]_{G_\ell+2S^2_\ell}\ar[ul]_{G_\ell+2S^1_\ell}\ar[ur]_{G_\ell+2S^3_\ell} & 
  \\  
 & *++\txt{$\Uone$}\ar[u]_{A_\ell} &
}
\]
\caption{Chain of subgroups $\SO(n)\subset\SO(9,1)\subset\SL(3,\OO)$.}
\label{Type123}
\end{center}
\end{minipage}
\end{center}
}
\end{figure}
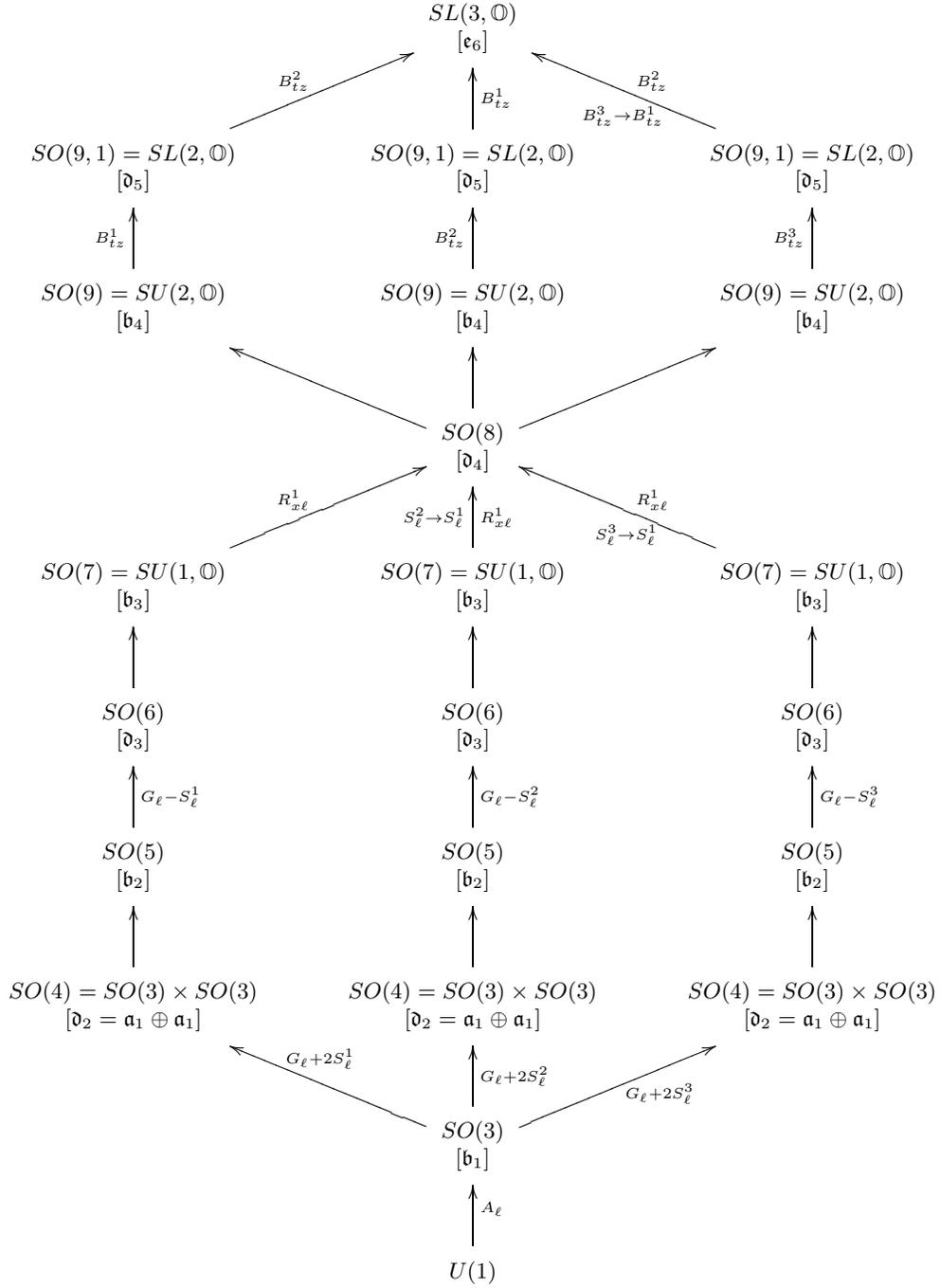

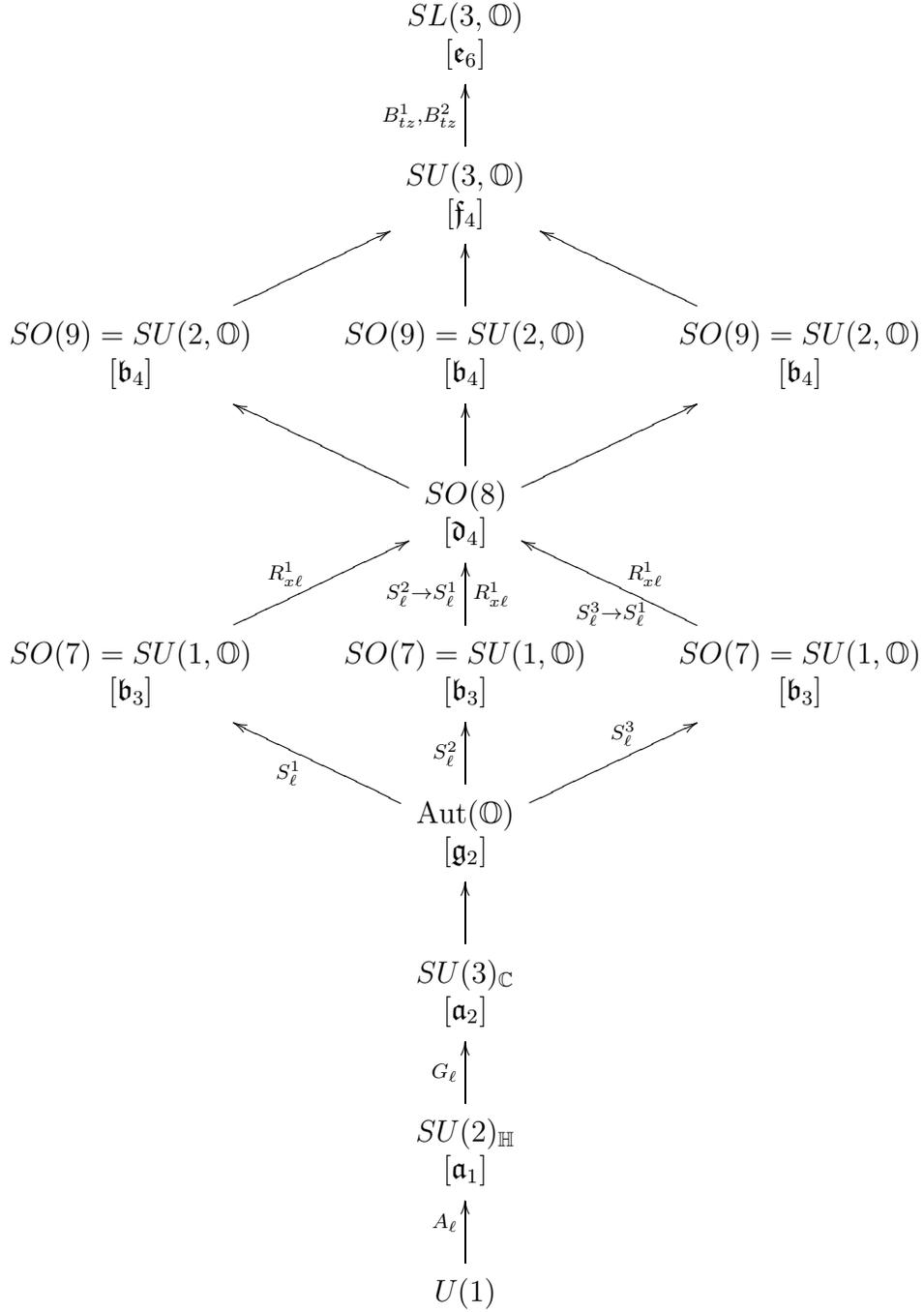
\begin{figure}[tbp]
\begin{center}
\begin{minipage}{6in}
\[
\xymatrix{
 & *++\txt{$\SL(3,\OO)$\\$[\ee_6]$} & \\
 & *++\txt{$\SU(3,\OO)$\\$[\ff_4]$}\ar[u]^{B^1_{tz}, B^2_{tz}} \\
*++\txt{$\SO(9) = \SU(2,\OO)$\\
  $[\bb_4]$}\ar[ur] & *++\txt{$\SO(9) = \SU(2,\OO)$\\
  $[\bb_4]$}\ar[u] & *++\txt{$\SO(9) = \SU(2,\OO)$\\
  $[\bb_4]$}\ar[ul] \\
 & *++\txt{$\SO(8)$\\$[\dd_4]$}\ar[u] \ar[ur] \ar[ul] \\
*++\txt{$\SO(7) = \SU(1,\OO)$\\
  $[\bb_3]$}\ar[ur]^{R^1_{x\ell}} & *++\txt{$\SO(7) = \SU(1,\OO)$\\
  $[\bb_3]$}\ar[u]^{S^2_\ell \to S^1_\ell}_{R^1_{x\ell}} &
	*++\txt{$\SO(7) = \SU(1,\OO)$\\
  $[\bb_3]$}\ar[ul]^{S^3_\ell \to S^1_\ell}_{R^1_{x\ell}} \\
 & *++\txt{$\Aut(\OO)$\\
	$[\gtwo]$}\ar[u]^{S^2_\ell} \ar[ur]^{S^3_\ell} \ar[ul]^{S^1_\ell} \\
 & *++\txt{$\SUC$\\$[\aa_2]$}\ar[u] \\
 & *++\txt{$\SUH$\\$[\aa_1]$}\ar[u]^{G_\ell} \\
 & *++\txt{$\Uone$}\ar[u]^{A_\ell}
}
\]
\caption{Chain of subgroups
$\SU(3)^C\subset\Gtwo\subset\SO(8)\subset\FF_4\subset\EE_6$.}
\label{Type123g2}
\end{minipage}
\end{center}
\end{figure}

We can use our particular choice of basis for the Lie algebra~$\ee_6$ to
identify two separate~$\SO(n)$ subgroup structures within the Lie group~$E_6$.
Figure~\ref{Type123} shows the~$\SO(n)$ subgroup chain of~$\SO(9,1)$ of type~$1$
in~$\SL(3,\OO)$, while Figure~\ref{Type123g2} shows the three~$\SO(9)$ subgroup
chains of~$F_4$ within~$E_6$.  In both subgroup structures, there is only
one~$\SO(8)$.  While~$G_2 \subset \SO(7)$, it is not a subset of~$\SO(6)$ in
Figure~\ref{Type123}.  Hence, we omit~$G_2$ from Figure~\ref{Type123}, but
include it in Figure~\ref{Type123g2} since our preferred basis for~$\SO(7)$
includes a basis for~$G_2$.  The figures indicate which Cartan element is
added to a group when it is expanded to a larger group, as well as giving the
classification of the corresponding Lie algebra.

\subsection{Type transformations}
\label{Type}

The \textit{discrete type transformation}~(\ref{typetran}) induced by~$\TT$
cyclically permutes Lorentz transformations of types~1, 2, or~3.  We have
\begin{align}
\TT^3 &= \II \\
\TT^\dagger &= \TT^{-1}
\end{align}
where $\II$ is the $3\times3$ identity matrix, and $\TT\in \SL(3,\OO)$, since
\begin{equation}
\det(\TT\XX\TT^\dagger) = \det(\XX)
\end{equation}
for $\XX\in\jalg$.  Although $\TT$ is not one of our elementary group
transformations, there are numerous identities of the form
\begin{align}
\TT &= R^1_{xz}(-\pi) \circ R^2_{xz}(-\pi) \nonumber\\
\TT &= R^2_{xz}(\pi) \circ R^1_{xz}(\pi)
	\circ R^2_{xz}(\pi) \circ R^1_{xz}(\pi) \\
\TT &= R^1_{xz}(\pi) \circ R^3_{xz}(\pi)
	\circ R^2_{xz}(\pi) \circ R^1_{xz}(\pi) \nonumber
\end{align}
These expressions make clear that $\TT\in \SL(3,\OO)$.  Furthermore, each of
these expressions may be expanded into a (different) \textit{continuous type
transformation} $\TT(\alpha)\in \SL(3,\OO)$ by letting the single fixed angle
($\pi$ or~$-\pi$) become arbitrary.  The resulting transformations are
\textit{not} one-parameter subgroups of~$\SL(3,\OO)$, but they do connect
transformations of different types.  We are therefore led to explore subgroups
of $\SL(3,\OO)$ that contain these (real!)\ type transformations, although it
suffices to consider subgroups containing $\TT$ itself.

\subsection{Type-independent subgroups}
\label{TypeSubs}

We list here some important groups which contain type transformations.  The
standard representation of~$\SO(3,\RR)$ is the group
\begin{equation}
\SO(3,\RR)_s = \langle R^1_{xz}, R^2_{xz}, R^3_{xz} \rangle
\end{equation}
This group obviously contains~$\TT$, as does the standard representation
\begin{equation}
\SL(3,\RR)_s
  = \langle R^1_{xz}, R^2_{xz}, R^3_{xz},
	B^1_{tz}, B^2_{tz}, B^1_{tx}, B^2_{tx}, B^3_{tx} \rangle
\end{equation}
of $\SL(3,\RR)$.  Using~$\ell$ as our preferred complex unit, we have the
standard representations
\begin{equation}
\SU(3,\CC)_s
  = \langle R^1_{xz}, R^2_{xz}, R^3_{xz},
	R^1_{x\ell}, R^2_{x\ell}, R^1_{z\ell}, R^2_{z\ell}, R^3_{z\ell} \rangle
\end{equation}
of~$\SU(3,\CC)$, and
\begin{equation}
\SL(3,\CC)_s
  = \SU(3,\CC)_s \cup \langle B^1_{tz}, B^2_{tz}, B^1_{tx},
	B^2_{tx}, B^3_{tx}, B^1_{t\ell}, B^2_{t\ell}, B^3_{t\ell} \rangle
\end{equation}
of $\SL(3,\CC)$.  These four groups are important because they contain the type
transformation~$\TT$.  If, for instance, some type~1 transformation~$R^1$ is
in a group~$G$ that has one of these groups as a subgroup, then $G$ must also
contain the corresponding type~2 and~3 transformations~$R^2$ and~$R^3$; we say
that $G$ is \textit{type independent}.

The standard representations~$\SO(3,\RR)_s$ and~$\SU(3,\CC)_s$ differ from our
preferred representations $\SOH=\SUH$ and~$\SUC$, which are subgroups
of~$G_2$.  For instance, the groups $\SU(3,\CC)_s$ and~$\SUC$ are both type
independent, but in~$\SU(3,\CC)_s$ the transformations~$R^1, R^2$ and~$R^3$ are
distinct while in~$\SUC$ the three transformations are equal; $\SUC$ does not
contain~$\TT$, nor does it need to.

\begin{figure}[tbp]
\[
\xymatrixcolsep{1pt}
\xymatrixrowsep{30pt}
\xymatrix@M=0pt{
*++\txt{$\sla(2,\OO)$ \\
  $\dd_5 $} & \txt{$\oplus$ \\
  $\oplus$} &
    *++\txt{$ \langle \dot B^2_{tz} - \dot B^3_{tz} \rangle $ \\
  $\dd_1$} \ar[d] \\
*++\txt{$\su(2,\OO)$ \\
  $\bb_4 $} \ar[u]  & \txt{$\oplus$ \\
  $\oplus$} &
    *++\txt{$ \langle \dot B^2_{tz} - \dot B^3_{tz} \rangle $\\
  $\dd_1$} \ar[d] \\
*++\txt{$\su(1,\OO)$ \\
  $\bb_3$} \ar[u] & \txt{$\oplus$ \\
  $\oplus$} &
    *++\txt{$ \langle \dot B^1_{tz}, \dot B^2_{tz} - \dot B^3_{tz} \rangle$ \\
  $ \dd_1 \oplus \dd_1$}\ar[d] \\
*++\txt{$\gtwo $ \\
  $\gtwo $} \ar[u] & \txt{$\oplus$ \\
  $\oplus$} & *++\txt{$ \sla(3,\RR)_s$ \\
  $\aa_1 \oplus \aa_1$} \ar[d] \\
*++\txt{$\suC $ \\
  $ \aa_2$} \ar[u]  & \txt{$\oplus$\\
  $\oplus$}  & *++\txt{$ \sla(3,\CC)_s$ \\
  $\aa_2 \oplus \aa_2$} \ar[d] \\
*++\txt{$\soH $ \\
  $\aa_1 $} \ar[u] & \txt{$\oplus$
  \\ $\oplus$} & *++\txt{$ \sla(3,\HH)$ \\
  $ \aa_5$} \\
}
\]
\caption{Type-dependent and type-independent subalgebras of~$\ee_6$.}
\label{TypeSub}
\end{figure}
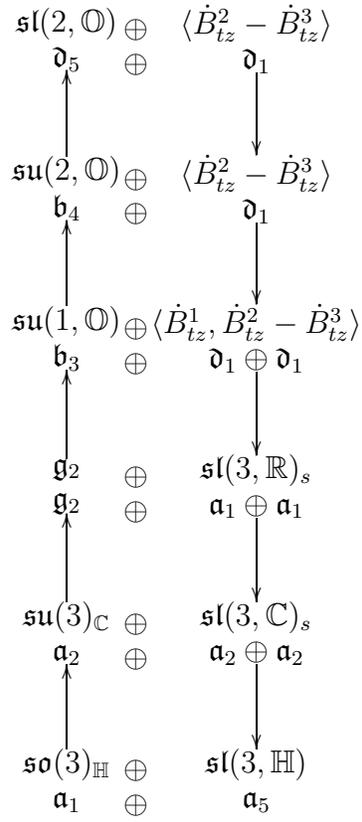

We use the type transformation~$\TT$ to provide insight into the structure of
the Lie algebra~$\sla(3,\OO)$.  The algebras~$\ga$ in the left column of
Figure~\ref{TypeSub} are subalgebras of the type~$1$ copy of~$\sla(2,\OO)$,
while each algebra~$\ga'$ in the right column is the largest subalgebra
of~$\sla(3,\OO)$ such that~$\ga\oplus\ga'$ is still simple.  When we
restrict~$\ga$ to a smaller subalgebra of~$\sla(2,\OO)$, it is sometimes
possible to expand the type-independent subalgebra~$\ga'$ to a larger
subalgebra of~$\sla(3,\OO)$.  Each arrow in the diagram indicates inclusion,
and a similar diagram holds for the corresponding subgroups of $\SL(3,\OO)$.

\subsection{\boldmath Reduction of~$\OO$ to~$\HH$,~$\CC$, and~$\RR$}
\label{Reduction}

We can also find subalgebras of~$\ee_6$ by restricting our generators to be
quaternionic, complex, or real.

Our preferred quaternionic subalgebra of $\OO$ is $\HH=\langle 1,k,k\ell,\ell
\rangle$, so we discard transformations involving $i$, $j$, $i\ell$, or
$j\ell$.  We therefore discard $3\times4=12$ boosts, $3\times4=12$ simple
rotations involving $z$, and $4$ simple rotations involving $x$ --- but we
must add back in 3 rotations involving $x$ of type~2, since we can no longer
use the middle relation in~(\ref{Srels}) to eliminate them.  Turning to the
transverse rotations, we need only consider transformations of type~1, and, as
discussed in Section~\ref{LorentzSub} (or after studying Table~\ref{AGS}), we
see that we must retain only the combinations $\dot G_q - \dot S^1_q$ for
$q\in\Im\HH$, thus discarding the remaining $21-3=18$ elements of $\so(7)$.  We
are left with $52-34+3=21$ rotations, and $26-12=14$ boosts.

The~$21$ compact generators form the algebra~$\su(3,\HH)$, a real form
of~$\cc_3=\sp(6)$, while all~$35$ together form~$\sla(3,\HH)$, a real form
of~\hbox{$\aa_5=\su(6,\CC)$}.  Restricting only to type~$1$ transformations,
we obtain~$10$ rotations and~$5$ boosts, thus reducing~$\sla(3,\HH)$
to~\hbox{$\sla(2,\HH)=\so(5,1)$}, a real form of~$\dd_3=\so(6)$, and
$\su(3,\HH)$ to~$\su(2,\HH)=\so(5)$, a real form of~\hbox{$\cc_2=\sp(4)$}.

Furthermore, the subalgebra
$\soH=\langle\dot A_k,\dot A_{k\ell},\dot A_\ell\rangle$ fixes $\HH$.  Thus,
for each of the above subalgebras $\ga$, we have $\ga\oplus\soH\in\sla(3,\OO)$.
In particular, \hbox{$\sla(3,\HH) \oplus \soH$} is therefore a subalgebra of
$\sla(3,\OO)$.

When restricting $\OO$ to our preferred complex subalgebra
\hbox{$\CC=\langle1,\ell\rangle$}, we obtain the classical Lie
algebras~$\su(3,\CC)_s$ and~$\sla(3,\CC)_s$ as previously discussed.  As there
is only one octonionic unit used to form~$\CC$, we do not need to use any of
the transformations from~$\SO(7)$, so we have 8 rotations and 8 boosts.  Using
all~$16$ transformations gives~$\sla(3,\CC)_s$, a real form
of~$\aa_2\oplus\aa_2=\su(3,\CC)\oplus \su(3,\CC)$ with~$8$ boosts, whereas we
obtain~$\su(3,\CC)_s$ by using only the~$8$ compact generators.  Further
restricting to the type~$1$ transformations reduces these two algebras
to~$\sla(2,\CC)_s=\so(3,1)_s$ and~$\su(3,\CC)_s$, which are real forms
of~\hbox{$\dd_2=\su(2,\CC)\oplus \su(2,\CC)$} and~$\aa_1=\su(3,\CC)$.

When we restrict~$\sla(3,\CC)_s$ to~$\sla(2,\CC)\subset \sla(2,\OO)$ (of type~1,
say), the smaller algebra no longer contains the type transformation~$\TT$,
but it does involve the octonionic direction~$\ell$.  Thus,
$\sla(2,\CC)\oplus \so(6)$, where~$\so(6)\subset \so(7)$ fixes~$\ell$, is also a
subalgebra of~$\sla(2,\OO)\subset \sla(3,\OO)$.

Finally, by restricting to real transformations, we are left with 3 rotations
and 5 boosts, which is a real form of~$\aa_2=\su(3,\CC)$ with~$5$ non-compact
elements.  This algebra may be further restricted to either~$\so(3,\RR)_s$,
whose group contains the type transformation, or~$\so(2,1)_s$, which is a
type~$1$ non-compact form of~$\aa_1=\so(3,\RR)$.

The above discussion of the result of restricting~$\sla(3,\OO)$
to~$\sla(n,\KK)$ for~$n = 1,2,3$ and \hbox{$\KK=\RR,\CC,\HH,\OO$} is
summarized in Figure~\ref{SubAlg}.  For each algebra~$\ga$ in
Figure~\ref{SubAlg}, we then list in Figure~\ref{PerpAlg} the maximal
subalgebra~$\ga'$ of~$\ee_6$ such that~$\ga\oplus\ga'\in\sla(3,\OO)$.
Here,~$\so(6)$ again denotes the subalgebra of type 1 which permutes
$\lbrace i,j,k,k\ell,j\ell,i\ell \rbrace$ but fixes~$\ell$.
Although~$\so(6)\not\subset\gtwo$, we do have~\hbox{$\suC\subset\so(6)$}.
We also write $\uonem$ for the non-compact real representation of $\dd_1$
generated by $\dotBtz$, which is discussed further in the next section.
%
Again, similar diagrams can be drawn for the corresponding subgroups of
$SL(3,\OO)$.

\begin{figure}[tbp]
\null\hspace{0.5in}
\begin{minipage}{5in}
\begin{center}
\xymatrixcolsep{8pt}
\xymatrix@M=1pt@H=0pt{
   & & & *++\txt{$\sla(3,\OO)$ \\ $[\ee_6]$} & & & \\
   & & & *++\txt{$\sla(3,\HH)$ \\ $[\aa_5]$}\ar[u] & & & \\
   & & & *++\txt{$\sla(3,\CC)_s$ \\ $[\aa_2\oplus \aa_2]$}\ar[u] & & & \\
  *++\txt{$\sla(2,\OO)$ \\ $[\dd_5]$} \ar[uuurrr] &
  *++\txt{$\sla(2,\HH)$ \\ $[\aa_3 = \dd_3]$} \ar[uurr] \ar[l] &
  *++\txt{$\sla(2,\CC)$ \\ $[\aa_1 \oplus \aa_1]$} \ar[ur] \ar[l] & &
  *++\txt{$\su(3,\CC)_s$ \\ $[\aa_2]$} \ar[ul] \ar[r] &
  *++\txt{$\su(3,\HH)$ \\ $[\cc_3]$} \ar[uull] \ar[r] &
  *++\txt{$\su(3,\OO)$ \\ $[\ff_4]$} \ar[uuulll] \\
   & & & *++\txt{$\su(2,\CC)$ \\ $[\aa_1]$} \ar[d] \ar[ul] \ar[ur] & & & \\
   & & & *++\txt{$\su(2,\HH)$ \\ $[\bb_2 = \cc_2]$}\ar[d] \ar[uurr] \ar[uull] & & & \\
   & & & *++\txt{$\su(2,\OO)$ \\ $[\bb_4]$}\ar[uuurrr] \ar[uuulll] & & & \\
}
\caption{Subalgebras~$\sla(n,\KK)$ and~$\su(n,\KK)$ of~$\sla(3,\OO)$.}
\label{SubAlg}
\end{center}
\end{minipage}
\end{figure}

\begin{landscape}
\begin{figure}[tbp]
\null\hspace{-0.25in}
{\Small
\xymatrixcolsep{12pt}
\xymatrixrowsep{12pt}
\xymatrix@M=2pt@H=2pt{
   & & & *++\txt{$\sla(3,\OO)\oplus 0 $ \\ $[\ee_6 \oplus 0]$}
	& & & \\
   & & & *++\txt{$\sla(3,\HH) \oplus \soH$ \\ $[\aa_5 \oplus \aa_1]$}
	\ar[u]<1ex>& & & \\
   & & & *++\txt{$\sla(3,\CC)_s \oplus \su(3)^C$ \\
	$[(\aa_2\oplus \aa_2) \oplus \aa_2 ]$}
	\ar[u]<1ex>& & & \\
  *++\txt{$\sla(2,\OO) \oplus \uonem$ \\ $[\dd_5 \oplus \dd_1]$}
	\ar[uuurrr]<1ex>&
  *++\txt{$\sla(2,\HH) \oplus \soH$ \\ $[(\aa_3 = \dd_3) \oplus \aa_1]$}
	\ar[uurr]<1ex> \ar[l]<-1ex> &
  *++\txt{$\sla(2,\CC) \oplus \so(6)$ \\ $[(\aa_1 \oplus \aa_1) \oplus \dd_3 ]$}
	\ar[ur]<1ex> \ar[l]<-1ex> & &
  *++\txt{$\su(3,\CC)_s \oplus \su(3)^C$ \\ $[\aa_2 \oplus \aa_2]$}
	\ar[ul]<1ex> \ar[r]<1ex>&
  *++\txt{$\su(3,\HH) \oplus \soH$ \\ $[\cc_3\oplus \aa_1]$}
	\ar[uull]<1ex> \ar[r]<1ex>&
  *++\txt{$\su(3,\OO) \oplus \uonem$ \\ $[\ff_4\oplus \dd_1]$}
	\ar[uuulll]<1ex>\\
   & & & *++\txt{$\su(2,\CC) \oplus \so(6)$ \\ $[\aa_1 \oplus \dd_3]$}
	\ar[d]<-1ex> \ar[ul]<1ex> \ar[ur]<1ex> & & & \\
   & & & *++\txt{$\su(2,\HH) \oplus \soH$ \\ $[(\bb_2 = \cc_2) \oplus \aa_1]$}
	\ar[d]<-1ex> \ar[uurr]<1ex> \ar[uull]<1ex> & & & \\
   & & & *++\txt{$\su(2,\OO)\oplus \uonem$ \\ $[\bb_4 \oplus \dd_1]$}
	\ar[uuurrr]<1ex>
	\ar[uuulll]<1ex> & & & \\
}
}
\caption{Subalgebras~$\sla(n,\KK) \oplus \ga'$ and~$\su(n,\KK)
\oplus \ga'$ of ~$\sla(n,\OO)$.}
\label{PerpAlg}
\end{figure}
\end{landscape}

\subsection{Subalgebras fixing type}
\label{FixEll}

Having just considered the subalgebras of $\gtwo$, and hence of $\ee_6$, that
leave invariant a preferred complex or quaternionic subalgebra of $\OO$, we
now ask what subalgebra of $\ee_6$ fixes all type~1 elements, that is, which
transformations leave $\matX$ alone in the first decomposition of
\hbox{$\XX\in\jalg$} shown in Table~\ref{Types}.  This subalgebra, which we
will call $\stab(I)$, turns out to be quite different from any of the others
discussed previously.

Clearly, no transformation in (type~1) $\sla(2,\OO)$ will be in $\stab(I)$.  We
therefore seek transformations of types~2 and~3.  Direct computation shows
that certain \textit{null rotations} will do the job.  Each of the 6 vector
spaces defined by
\begin{equation}
b^a_\pm = \langle \dot B^a_{tx} \mp \dot R^a_{xz},
		\dot B^a_{tq} \pm \dot R^a_{zq} \rangle
\end{equation}
is in fact an abelian subalgebra of $\sla(3,\OO)$, and in each case the given
basis elements are null according to the Killing form --- the Killing
form is in fact identically zero on each of these subalgebras.  Each of these
subalgebras fixes all elements of a particular type; we have
\begin{equation}
\stab(I) = b^2_+ \oplus b^3_-
\end{equation}
with cyclic permutations holding for $\stab(II)$ and $\stab(III)$.

Since $\stab(I)$ contains no elements of (type~1) $\sla(2,\OO)$, we expect that
$\sla(2,\OO)\oplus\stab(I)$ will be a $45+16=61$-dimensional subalgebra of
$\sla(3,\OO)$.  Checking commutators, this turns out to be correct, but with an
unexpected surprise: $\stab(I)$ is an ideal of $\sla(2,\OO)\oplus\stab(I)$, so
this subalgebra is neither simple nor semisimple.

If we further define
\begin{equation}
\stab(I)^\perp = b^2_- \oplus b^3_+
\end{equation}
to be the 16 null rotations of types~2 and~3 that are not in $\stab(I)$, then
we have the intriguing decomposition 
\begin{equation}
\sla(3,\OO) = \sla(2,\OO) \oplus \stab(I) \oplus \stab(I)^\perp \oplus \uonem
\end{equation}
with $\uonem$ again denoting the non-compact real representation of $\dd_1$
generated by $\dotBtz$.

We can now easily determine the subalgebras of $\ee_6$ that, say, leave $\HH$
or $\CC$ in type 1 elements invariant.  All we have to do is combine the
relevant subalgebra of $\sla(2,\OO)$ --- in this case $\suH$ or $\suC$,
respectively --- with $\stab(I)$.  Each such algebra, here
$\suH\oplus\stab(I)$ and $\suC\oplus\stab(I)$ is a subalgebra of $\ee_6$
which, however, is neither simple nor semisimple.  Two further examples are
the 52-dimensional subalgebras $\su(2,\OO)\oplus\stab(I)$, which fixes
(type~1)~$t$, and $\so(8,1)_\ell\oplus\stab(I)$, where $\so(8,1)_\ell$ fixes
(type~1)~$\ell$ (and therefore does not contain $g_2$).

\newpage

\section{Conclusion}
\label{conclusion}

In this paper, we have given an explicit description of the subgroup structure
of $\SL(3,\OO)$, based on the ``type'' structure inherent in the embedding of
$\SL(2,\OO)$ in $\SL(3,\OO)$, and on the structure of $\SL(2,\OO)$ itself.  In
the process, we have provided explicit realizations of some of the remarkable
properties of $G_2$.  The internal structure of $G_2$, such as the $\SU(3)$ and
$\SU(2)$ subgroups fixing either a complex or quaternionic subalgebra, may be
especially relevant to attempts to use $\SL(3,\OO)$ to describe fundamental
particles, as discussed further in~\cite{York}.  Furthermore, we have seen
explicitly how $G_2$ is preserved under triality, as discussed
in~\cite{Denver}.  Finally, we have constructed the groups leaving the type
structure invariant, which we suspect may play a prominent role in describing
the interactions of fundamental particles.

However, the story is only partially complete.  There are other interesting
subgroups of $\SL(3,\OO)$, closely related to the 4 other real forms of $E_6$.
In particular, we have not yet identified any of the $C_4$ subgroups of $E_6$.
In other work~\cite{E6cartan}, we extend, and in a sense complete, the present
investigation by constructing and discussing chains of subgroups adapted to
these other subgroups.  We hope that the resulting maps of $E_{6(-26)}$ will
prove useful in further attempts to apply the exceptional groups to nature.


\section*{Acknowledgments}

This paper is a revised version of Chapter~4 of a dissertation submitted by AW
in partial fulfillment of the degree requirements for his Ph.D.\ in
Mathematics at Oregon State University~\cite{aaron_thesis}.  The revision was
made possible in part through the support of a grant from the John Templeton
Foundation.

\newpage

\bibliographystyle{unsrt}
\bibliography{e6sub}

\end{document}